\documentclass[letterpaper, 10 pt, conference]{ieeeconf}  % Comment this line out if you need a4paper

\IEEEoverridecommandlockouts                              % This command is only needed if
\usepackage{amsmath} % assumes amsmath package installed
\usepackage{amssymb}  % assumes amsmath package installed
\usepackage{xcolor}
\usepackage{graphicx}
\usepackage[ruled]{algorithm2e}
\usepackage{cite}

\newcommand{\mb}[1]{\mathbf{#1}}
\newcommand{\x}{\mb{x}}
\newcommand{\tx}{\tilde{\x}}
\newcommand{\xdot}{\dot{\mb{x}}}
\newcommand{\dx}{\dot{x}}
\newcommand{\X}{\mb{X}}
\newcommand{\tX}{\tilde{\X}}
\newcommand{\m}{\mb{m}}
\newcommand{\y}{\mb{y}}
\newcommand{\ydot}{\dot{\y}}

\newcommand{\q}{\mb{a}}
\newcommand{\f}{\mb{f}}

\newcommand{\E}[1]{\mb{E}[#1]}
\newcommand{\Dx}{\mb{\mathcal{X}}}
\newcommand{\nx}{d}
\newcommand{\nb}{n_B}
\newcommand{\nd}{n_D}
\newcommand{\Real}{\mathbb{R}}
\newcommand{\maxent}{\textit{maxent }}
\newcommand{\cvxh}[1]{\text{Conv}(#1)}
\newcommand{\fhat}{\hat{f}}

\newcommand{\eqnlabel}[1]{\label{eqn:#1}}
\newcommand{\figlabel}[1]{\label{fig:#1}}
\newcommand{\eqn}[1]{eq.(\ref{eqn:#1})}

\newcommand{\fig}[1]{Fig.\ref{fig:#1}}
\newcommand{\Fig}[1]{Fig.\ref{fig:#1}}

\title{\LARGE \bf
Surrogate Modeling of Dynamics From Sparse Data Using Maximum Entropy Basis Functions
}

\author{Vedang M. Deshpande$^{1}$ and Raktim Bhattacharya$^{2}$% <-this % stops a space
\thanks{This work was supported by the National Science Foundation (grant number: 1762825).}% <-this % stops a space
\thanks{$^{1}$Vedang M. Deshpande is a PhD student in Aerospace Engineering, Texas A\&M University, College Station, TX 77843, USA.
        {\tt\small vedang.deshpande@tamu.edu}}%
\thanks{$^{2}$Raktim Bhattacharya is Associate Professor in Aerospace Engineering, Electrical \& Computer Engineering, Texas A\&M University, College Station, TX 77843, USA.
        {\tt\small raktim@tamu.edu}}%
}

\begin{document}

\maketitle
\thispagestyle{empty}
\pagestyle{empty}

%%%%%%%%%%%%%%%%%%%%%%%%%%%%%%%%%%%%%%%%%%%%%%%%%%%%%%%%%%%%%%%%%%%%%%%%%%%%%%%%
\begin{abstract}
In this paper we present a data driven approach for approximating dynamical systems. A dynamics is approximated using basis functions, which are derived from maximization of the information-theoretic entropy, and can be generated directly from the data provided. This approach has advantages over other methods, where a dictionary of basis functions have to be provided by the user, which is non trivial in some applications. We compare the accuracy of the  proposed data-driven modeling approach to existing methods in the literature, and demonstrate that for some applications the maximum entropy basis functions provide significantly more accurate models.
\end{abstract}

%%%%%%%%%%%%%%%%%%%%%%%%%%%%%%%%%%%%%%%%%%%%%%%%%%%%%%%%%%%%%%%%%%%%%%%%%%%%%%%%
\section{INTRODUCTION}
It is often difficult to model complex systems from fundamental governing laws, as such models typically have high computational complexity, if at all possible. Some of these systems are often built before any formal modeling effort, and data is available from the various experiments performed on these systems. It is also possible that it is expensive to run experiments on these systems, or may be impossible to get data for certain operating conditions. Consequently, the available data is sparse. For example, in diverse applications such as aeroelastic tailoring of wind turbine blades, shock interactions in scramjet engine dynamics, and mechano-biological response of cells, modeling from first principles is a daunting task as there are many unknowns. However, it is possible to obtain data from limited number of experiments. Therefore, determining surrogate models from sparse data will be extremely useful, and will be useful for various model-based analysis and design activities, such as uncertainty quantification and control system design.

Polynomial interpolation is one of many techniques used for constructing approximate models based upon given data. Depending upon the problem and available data, different types of polynomials (e.g. Lagrange, Chebyshev) are used as basis functions for constructing interpolants \cite{polyBook}.

One class of interpolation problems, known as polygonal interpolation, deals with the construction of interpolant within a convex polytope if the function values at the vertices are known. Each node is associated with a shape function, which is used in the construction of the interpolant. The choice of shape function depends upon the desired properties of the interpolant.

Maximum entropy based shape functions were introduced for polygonal interpolation in \cite{sukumar}.
The principle of maximum entropy (\textit{maxent}) which is discussed in the subsequent section, was introduced in \cite{jaynes1, jaynes2} as a method to draw the least biased inference from limited or insufficient data.
Extensions of such maximum entropy based interpolation schemes were studied in various works such as \cite{ortiz, sukumarRelEntr1, sukumarRelEntr3, sukumarRelEntr2}. In this paper, we explore maximum entropy functions as bases for constructing an approximant of the unknown function or dynamics.

Polygonal interpolation schemes use scattered data points (basis nodes) and corresponding function values to construct a convex approximant. Each basis node is associated with a shape function, which is evaluated for a point at which approximation is to be made. Therefore, the point of interest where one wish to evaluate the approximant must lie within the convex hull of the available nodes. In the present work, we relax this constraint by choosing a set of basis nodes independently of the available data points, and given function values are transformed w.r.t. basis nodes by minimizing a suitably defined cost function. Approximants constructed using this method can be evaluated for a point, which may not lie within the convex hull of available data points. However, needless to say, the point must lie within the convex hull of the independently selected basis nodes.

Popular data driven modeling approaches, such as \cite{sindy}, utilize a user defined library of bases (e.g. polynomials, trigonometric functions)  and minimize $l_2$ error for given data to determine coefficients associated with the basis functions. If the available data lies within the span of specified bases, then the $l_2$ convergence is achieved. Otherwise, as we show in examples in subsequent sections, approximation errors can be significant.
Moreover, it is well known that higher order polynomial interpolations suffer from over fitting of the given data, also known as Runge's phenomenon. The approach presented in this paper utilizes \maxent functions as bases which are derived from the given data for function approximation. The approach is extended for dynamics modeling using sparse observations.

The key contributions of this work are highlighted below.
\begin{enumerate}
\item[$\bullet$] This work extends the classical polygonal interpolation schemes by allowing basis nodes to be selected independently of the available data points.
\item[$\bullet$] The framework does not rely on the user defined library of basis functions. They are derived from the available data set and selected basis nodes by maximization of the entropy. This is particularly useful when no information is available about the equations which govern the given data.
\item[$\bullet$] Through multiple examples, we demonstrate the application of this framework for various systems.  We show that this framework provides quite accurate approximants even if the available data set is sparse.
\end{enumerate}

The organization of the paper is as follows. We begin Section \ref{maxent_fcn} with the principle of maximum entropy (\textit{maxent}) followed by a brief discussion on polygonal interpolation, and global and local \maxent basis functions or barycentric coordinates.
Section \ref{fcn_approx} builds upon the preliminaries from Section \ref{maxent_fcn} to develop a framework for function approximation and few examples are discussed. Section \ref{dyn_model} presents an algorithm which summarizes the function approximation framework in the context of dynamics modeling. Concluding remarks are presented in Section \ref{conclude}.
\section{MAXIMUM ENTROPY BASIS FUNCTIONS} \label{maxent_fcn}

\subsection{Principle of maximum entropy}
Let $\mb{p} := [p_1, p_2 \cdots p_N]^T$ be the discrete probability distribution associated with mutually independent events $x_i$ for $i=1\cdots N$. Therefore, $p_i \geq 0$ and $\sum_{i=1}^N p_i = 1$. Then the information-theoretic entropy of the distribution $\mb{p}$ is defined as \cite{jaynes1,jaynes2}
\begin{align}
  H(\mb{p}) := -\sum_{i=1}^N p_i \log(p_i), \eqnlabel{entropy}
\end{align}
where $p_i \log(p_i):=0$ if $p_i=0$.

Suppose the probability distribution $\mb{p}$ is unknown. However, expected value of a function $\E{\mb{g}(x)}:=\sum_{i=1}^N p_i \mb{g}(x_i)$ is given. The principle of maximum entropy states that, out of all possible distributions which represent the given data, the least biased or the probability distribution ($\mb{p}^*$) with the maximum likelihood, is the one which maximizes the entropy defined in \eqn{entropy}, i.e.,
\begin{subequations}
\begin{equation}
  \mb{p}^* = \operatorname*{arg\,max}_{\mb{p}} H(\mb{p}),
\end{equation}
\begin{align}
  \text{such that } \sum_{i=1}^N p_i \mb{g}(x_i) &= \E{\mb{g}(x)} , \text{ and } \sum_{i=1}^N p_i = 1.
\end{align}
  \eqnlabel{constraintME}
\end{subequations}
Reference \cite{sukumar} drew an analogy between  the problem of polygonal interpolation and the maximum entropy probability distribution. We briefly discuss the concept of polygonal interpolation below.

\subsection{Maximum entropy and polygonal interpolation}
Let $S:=\{\x_i\}_{i=1}^N \subset \Dx \subset \Real ^{\nx}$, be a set of $N$ discrete points in $\nx-$dimensional space and $\Dx$ be compact. Then the convex hull of the set $S$ is defined as
\begin{align*}
\text{Conv}(S):= \{ \x | \x = \sum_{i=1}^N w_i \x_i, \sum_{i=1}^N w_i =1 , w_i \geq 0, \x_i \in S \} .
\end{align*}
Each point or node $\x_i$ is associated with a basis or shape function $\psi_i(\cdot) \geq 0$. Let $f(\x): \Dx \rightarrow \Real$ be a scalar function. For any $\x \in \text{Conv}(S)$, the polygonal interpolant $\hat{f}(\cdot)$ of $f(\cdot)$ at $\x$ is given by
\begin{align}
  \hat{f}(\x) = \sum_{i=1}^N \psi_i(\x) f(\x_i). \eqnlabel{polyInterp}
\end{align}
Coefficients $\psi_i(\x)$ are also known as barycentric coordinates of $\x$ w.r.t. nodes $\x_i \in S$. Constraints are imposed on $\psi_i$ so that $\hat{f}(\cdot)$ can exactly recover constant and linear functions. These constraints are given by
\begin{subequations}
\begin{equation}
  \sum_{i=1}^N \psi_i(\x) = 1,
\end{equation}
\begin{equation}
  \sum_{i=1}^N \psi_i(\x) \x_i = \x, \text{ or, } \X \mb{\Psi}(\x) = \x,
  \eqnlabel{linearConstr}
\end{equation}
\eqnlabel{linearConstrAll}
\end{subequations}
where $\X := [\x_1, \x_2 \cdots \x_N]$, and $\mb{\Psi}(\x):= [\psi_1(\x), \psi_2(\x) \cdots \psi_N(\x)]^T$.

Comparison of these equations with \eqn{constraintME} yields that $\psi_i(\x)$ can be interpreted as probability associated with the \textit{event} $\x_i$, and here $\mb{g}(\x) = \x$. Therefore, $\mb{\Psi}(\x)$ which maximize the entropy can be interpreted as the least biased barycentric coordinates of $\x$ w.r.t. nodes $\x_i \in S$.
For numerical stability of the algorithm, it is customary to use shifted coordinate system with origin at $\x$. Therefore, we define $\tx_i:= \x_i - \x$ and the constraints in \eqn{linearConstrAll} can be rewritten as
$\tX \mb{\Psi}(\x) = \mb{0}$, where $\tX:=[\tx_1, \tx_2 \cdots \tx_N]$.
The entropy maximization problem is formally presented as follows.
\begin{subequations}
\begin{equation}
  \mb{\Psi}(\x)^* = \operatorname*{arg\,max}_{\mb{\Psi}(\x)} H(\mb{\Psi}(\x)),\\
\end{equation}
\begin{equation}
  \text{such that } \tX \mb{\Psi}(\x) = \mb{0} %\text{ and } \sum_{i=1}^N \mb{\Psi}(\x)_i = 1.
  \eqnlabel{PhiMEConstr}
\end{equation}
  \eqnlabel{PhiME}
\end{subequations}
This is a convex optimization problem with linear constraints. Solution to this problem using the method of Lagrange multipliers is discussed in \cite{jaynes1,sukumar}.
The $i^{\text{th}}$ basis function in terms of Lagrange multipliers is given by
\begin{align}
\psi_i(\x) = \frac{e^{-\boldsymbol{\lambda}^T \tx_i}}{\sum_{j=1}^N e^{-\boldsymbol{\lambda}^T \tx_j}}
\eqnlabel{psi}
\end{align}
where $\boldsymbol{\lambda}^T:=[\lambda_1, \lambda_2 \cdots \lambda_{\nx}]$ are Lagrange multipliers associated with equality constraints in \eqn{PhiMEConstr}.

Determining Lagrange multipliers reduces to solving following system of nonlinear equations in terms of $\lambda_i$.
\begin{align}
 -\sum_{i=1}^N \tx_i e^{-\boldsymbol{\lambda}^T \tx_i}  = \mb{0}. \eqnlabel{lamEqns}
\end{align}
Because of the nonlinearity, numerical solvers are employed to determine optimal Lagrange multipliers.
Corresponding basis functions $\mb{\Psi}(\x)^*$ are recovered using \eqn{psi}. These basis functions are called \textit{global maxent} basis functions.
Interpolants constructed using such global functions often lead to poor fit to the given data \cite{ortiz}. Therefore, it is desirable to have control over the degree of locality of basis functions. This is achieved by using \textit{local maxent} basis functions which is discussed next.

\subsection{Local maximum entropy basis functions}
The decay of a basis function with distance away from its corresponding node, also known as width of the basis function is controlled by modifying the cost function in optimization problem given by \eqn{PhiME}. A generalized formulation to derive local \textit{maxent} basis functions is discussed in detail in \cite{ortiz}.
For the purpose present study, we adopt Gaussian radial basis functions to define a prior and minimize relative entropy to determine local \textit{maxent} basis functions \cite{sukumarRelEntr1,sukumarRelEntr2}.

For any $\x \in \text{Conv}(S)$, the Gaussian prior, $\m(\x)$, is defined as
\begin{align}
\m(\x):= [m_1, m_2 \cdots m_N]^T \text{ and } m_i(\x):= e^{-\beta ||\tx_i||_2^2},
\eqnlabel{prior}
\end{align}
where $\beta \geq 0$. Relative or cross entropy is defined as
\begin{align}
  \bar{H}\big(\mb{\Psi}(\x), \m(\x)\big) := \sum_{i=1}^N \psi_i \log\Big(\frac{\psi_i}{m_i}\Big). \eqnlabel{relEntropy}
\end{align}
We denote local \maxent basis functions by $\bar{\mb{\Psi}}(\x)$ and they are determined by solving the following optimization problem
\begin{subequations}
\begin{equation}
  \bar{\mb{\Psi}}(\x)^* = \operatorname*{arg\,min}_{\bar{\mb{\Psi}}(\x)} \bar{H}\big(\bar{\mb{\Psi}}(\x), \m(\x)\big),\\
\end{equation}
\begin{equation}
  \text{such that } \tX \bar{\mb{\Psi}}(\x) = \mb{0}.
  \eqnlabel{PhiMEConstr_loc}
\end{equation}
  \eqnlabel{PhiME_loc}
\end{subequations}
Similar to \eqn{lamEqns}, this optimization problem reduces to following system of nonlinear equations in Lagrange multipliers $\lambda_i$,
\begin{align}
 \sum_{i=1}^N \tx_i m_i(\x) e^{-\boldsymbol{\lambda}^T \tx_i}  = \mb{0}. \eqnlabel{lamEqns_loc}
\end{align}
Solution of \eqn{lamEqns_loc} is determined numerically. We use \texttt{fsolve} function in \texttt{MATLAB} for this purpose.
Local basis functions in terms of $\boldsymbol{\lambda}$ are recovered as
\begin{align}
\bar{\psi}_i(\x) = \frac{m_i(\x) e^{-\boldsymbol{\lambda}^T \tx_i}}{\sum_{j=1}^N m_j(\x) e^{-\boldsymbol{\lambda}^T \tx_j}}.
\eqnlabel{psi_loc}
\end{align}
Note, parameter $\beta$ controls the degree of locality of basis functions. Higher the value of $\beta$, larger is the decay of basis function away from its corresponding node.
Thus, locality of basis functions increases with the value of $\beta$. For uniform prior, i.e. $m_i(\x) = 1/N$, for $i = 1,2\cdots N$, or $\beta = 0$, we recover global \maxent basis functions given by \eqn{psi}.
Therefore, for the simplicity of notation, hereafter we use  $\mb{\Psi}(\x)$ or $\psi_i(\x)$ to denote the solution of \eqn{PhiME_loc}, i.e. optimal local \maxent basis functions evaluated at $\x$.

\section{FUNCTION APPROXIMATION USING LOCAL MAXENT BASIS FUNCTIONS} \label{fcn_approx}
In classical polygonal interpolation, basis functions, $\mb{\Psi}(\x)$, are evaluated for a given $\x \in \cvxh{S}$ by solving optimization problem defined by \eqn{PhiME_loc}. Then the interpolant of a function at $\x$ is given by \eqn{polyInterp} \cite{sukumarRelEntr1,sukumarRelEntr2}.
For one to be able to use such interpolant over the entire domain $\Dx$ of the function, $\Dx$ must be the convex hull of the given data set $S$, because optimization problem \eqn{PhiME_loc} is feasible only if $\x \in \cvxh{S}$.
This condition is not necessarily satisfied for real world systems, especially if the available data is sparse. To this end, we propose the following formulation for function approximation.

Let $B:=\{\x_i\}_{i=1}^{\nb} \subset \Dx \subset \Real ^{\nx}$, such that $\Dx = \cvxh{B}$. Therefore, for any $\x \in \Dx$, \eqn{PhiME_loc} is feasible and we can determine barycentric coordinates of $\x$ w.r.t. nodes in $B$.
Hereafter, unless specified, $\mb{\Psi}(\x)$ will be used to denote optimal local \maxent basis functions or barycentric coordinates of $\x$ w.r.t. $\x_i \in B$. Therefore, $\mb{\Psi}(\x)$ is an $\nb-$dimensional vector.

Without loss of generality, let $f(\x):\Dx \rightarrow \Real$ be a real valued scalar function.
Let $D:=\{\y_i\}_{i=1}^{\nd} \subset \Dx \subset \Real ^{\nx}$, be the set of points for which functions values $f(\y_i)$ are given. Set $D$ need not be same as set $B$. Based upon this information, our objective is to construct an approximant $\fhat(\cdot)$ of $f(\cdot)$.

Instead of interpolating given function values $f(\y_i)$, we use $\psi_i(\x)$ as basis functions to approximate $f(\x)$ and minimize $l_2$ error over the given data set. In particular,
\begin{align}
\fhat(\x) := \sum_{i=1}^{\nb} a_i^* \psi_i(\x) = \q^{*T}\mb{\Psi}(\x)
\eqnlabel{fhat}
\end{align}
where $[a_1^*, a_2^* \cdots a_{\nb}^*]^T=:\q^*$ is the solution of following optimization problem
\begin{align}
  \q^* = \operatorname*{arg\,min}_{\q} \Big(||\varepsilon||_2 + \alpha||\q||_1 \Big)
  \eqnlabel{opt_a}
\end{align}
where
\begin{align*}
 \alpha \geq 0,  \text{ and } \varepsilon := \begin{pmatrix}f(\y_1)\\ f(\y_2)\\ \vdots \\f(\y_{\nd}) \end{pmatrix}  -  \begin{pmatrix} \q^T\mb{\Psi}(\y_1)\\ \q^T\mb{\Psi}(\y_2) \\ \vdots \\ \q^T\mb{\Psi}(\y_{\nd}) \end{pmatrix}.
\end{align*}
Note that $\varepsilon$ is the error vector and $\fhat(\y_i) = \q^T\mb{\Psi}(\y_i)$ is the approximated value of $f(\y_i)$, where $\mb{\Psi}(\y_i)$ are the local \maxent basis functions evaluated for $\y_i$ w.r.t. nodes in $B$.
$l_1$ regularization on $\q$ prevents over fitting of the given data and $\q^*$ from taking very large values.
Optimization problem given by \eqn{opt_a} is convex and can be solved numerically using software packages such as \texttt{cvx} \cite{cvx}.
We demonstrate the application of this framework using following examples.
\subsection*{Examples}
\subsubsection{$f(x) = \sin(2\pi x)$} $x \in [0,1]$. We select uniformly spaced $\nb=10$ points in $[0,1]$ as nodes for \maxent basis functions.
$\nd = 20$ random data points sampled uniformly from [0,1] and their corresponding function values serve as our \textit{given} or \textit{training data set}. Using the formulation above, we construct an approximant $\fhat(x)$ and evaluate it for the \textit{test data set} which is uniformly spaced $50$ points in $[0,1]$.
In \fig{sinx}, $\fhat(\cdot)$ evaluated for test data is shown by red dashed line and true function values are shown by black dotted line. Similarly, red circles denote the $\fhat(\cdot)$ evaluated at the given data points and black crosses denote true function values at these points. Clearly, $\fhat(\cdot)$ approximates the given function $f(x) = \sin(2\pi x)$  with good accuracy.
Root mean square (RMS) error in the approximation for the given data set is $\mathcal{O}(10^{-4})$ and for the test data set is $\mathcal{O}(10^{-3})$. These results are obtained for $\beta = 100$ and $\alpha = 0$.
\begin{figure}[thpb]
   \centering
   %\framebox{\parbox{3in}{\includegraphics[width=0.5\textwidth]{sinx.eps}}}
   \includegraphics[width=0.5\textwidth]{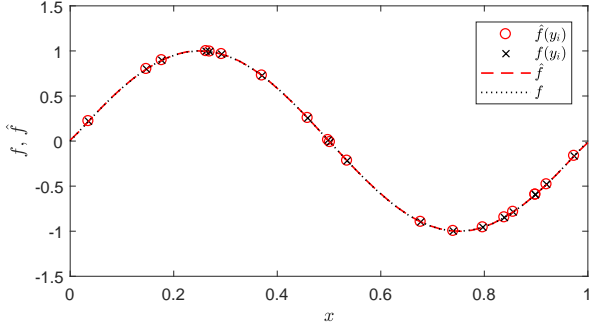}
   \caption{Comparison of approximant and true function for $f(x)=\sin(2\pi x)$.}
   \figlabel{sinx}
\end{figure}

\subsubsection{$f(\x) = 2x_1 e^{-4||\x||_2^2}$} $\x:=[x_1, x_2]^T \in [0,1] \times [0,1]$. \label{expx}
We use a uniformly spaced grid of $8\times 8$ ($\nb = 8^2$) points as nodes for basis functions. Another grid of $16 \times 16$ ($\nd = 16^2$) points and their corresponding function values are used as training data set.
We test the accuracy of approximated function on the grid of $32 \times 32$ points. \Fig{xexpx} shows the results obtained for $\beta = 10$ and $\alpha = 0$. RMS error for both training and test data points is $\mathcal{O}(10^{-4})$.
\begin{figure}[thpb]
   \centering
   \includegraphics[width=0.5\textwidth]{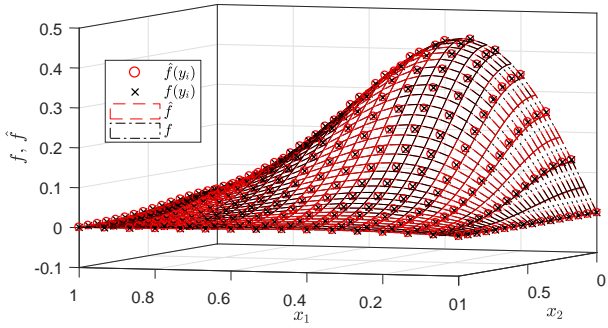}
   \caption{Comparison of approximant and true function for $f(\x) = x_1 e^{-||\x||_2^2}$.}
   \figlabel{xexpx}
\end{figure}

\subsubsection{$f(\x) = (1-x_1)^2 + 100(x_2-x_1^2)^2$} $\x:=[x_1, x_2]^T \in [-1,1] \times [-1,1]$. This is known as Rosenbrock function \cite{rosenbrock} which is often used as a test problem for optimization algorithms.
We use similar grids as the previous example constructed over $[-1,1] \times [-1,1]$, i.e. $\nb = 8^2, \nd = 16^2$ and approximated function is tested on the grid of $32^2$ points. RMS error for both training and test data is $\mathcal{O}(10^{-2})$.
These results  are obtained for $\beta =5$, $\alpha = 0$ and shown in \fig{rosen}.
\begin{figure}[thpb]
   \centering
   \includegraphics[width=0.5\textwidth]{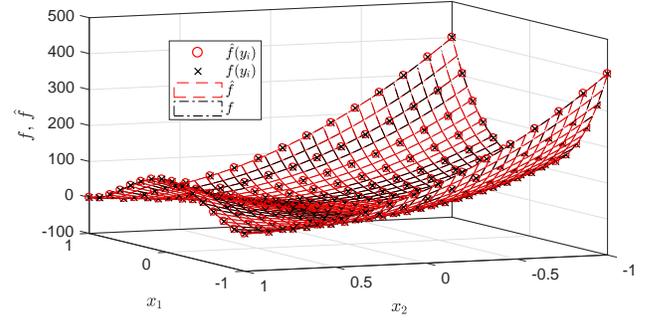}
   \caption{Comparison of approximant and true Rosenbrock function.}
   \figlabel{rosen}
\end{figure}

These examples show that local \maxent basis functions can approximate nonlinear functions with reasonable accuracy. We empirically observed that if given data points are used as basis nodes, i.e. $B = D$, and $\beta$ is sufficiently large, i.e. basis functions exhibit greater degree of locality, then the RMS error corresponding to given points is less than $10^{-14}$.
However, order of magnitude of RMS error for the test data set remains unchanged. As shown in examples, this allows us to use sparser set of nodes for basis. e.g. in example \ref{expx}, instead of using $16^2$ given data points, we use $8^2$ points as nodes for basis functions while maintaining the similar order of accuracy of the approximated function.
The information available from $16^2$ data points is transformed in reduced $8^2$-dimensional vector $\q^*$. This greatly reduces the number of basis functions, hence, simplifying the approximation given by \eqn{fhat}.

We compare our results with the algorithm SINDy \cite{sindy}. SINDy recovers the Rosenbrock and $\sin(\cdot)$ function with testing RMS error of $\mathcal{O}(10^{-14})$, as they lie in the span of its specified library of $4^{\text{th}}$ order polynomial and trigonometric basis functions.
However, error for the function in example \ref{expx} is an order of magnitude greater than the \maxent based approximant. Therefore, these preliminary results show that if no information is available about the nature of true function, i.e. if the given data does not lie within the span of user defined basis functions,
then \maxent based approximants can be employed for better accuracy.

\subsection*{Effect of $\alpha$}
In the examples discussed above, we have used $\alpha = 0$ which is the coefficient of $l_1$ regularization on $\q$ in optimization problem \eqn{opt_a}. Here, we briefly demonstrate the effect of $\alpha$ using $\sin(\cdot)$ function example.
If the available data points are clustered in the domain as shown in \fig{sinxalpha}, then approximant can suffer from large errors at the boundaries of the domain for $\alpha = 0$, represented by red dashed line. Blue dashed line corresponding to $\alpha = 10^{-3}$ shows that with $l_1$ regularization on $\q$, such errors can be reduced.

Moreover, $\alpha$ can be used as a parameter to induce sparsity in the set of basis nodes. For nonzero $\alpha$, if magnitude of some elements of $\q^*$ are very small relative to other elements, then it implies that the basis functions corresponding to those $a_i^*$ are not required. Thus, we can reduce the number of basis nodes to make the set of basis functions sparse. e.g. In \fig{sinxalpha}, red and blue lines correspond to $\nb=20$.
The green dashed line corresponds to sparser basis set with $\nb = 10$ and it actually improves the accuracy of the approximant. Therefore, in addition to $\beta$, $\alpha$ can also be used as a control parameter to improve the accuracy of the approximant.

\begin{figure}[thpb]
   \centering
   \includegraphics[width=0.5\textwidth]{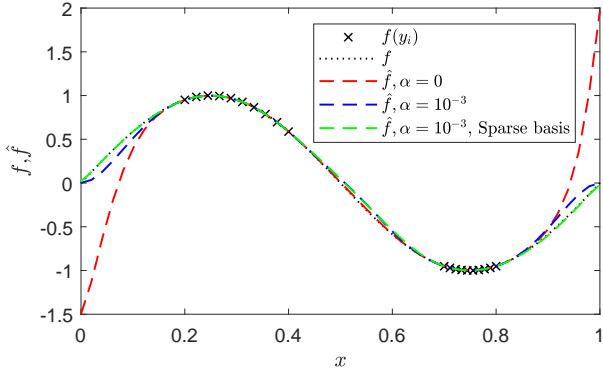}
   \caption{Comparison of approximant and true function for $f(x)=\sin(2\pi x)$ for different values of $\alpha$, and fixed $\beta = 100$.}
   \figlabel{sinxalpha}
\end{figure}

With these remarks on $\alpha$, we conclude our discussion on function approximation using \maxent basis functions. In the text to follow, we extend the current formulation for data driven modeling of dynamics.

\section{MODELING OF DYNAMICS} \label{dyn_model}
Let us consider a dynamical system governed by
\begin{align}
\xdot = \f(\x), \eqnlabel{dyn}
\end{align}
where, $\x :=[x_1, x_2 \cdots x_{\nx} ]^T \in \Dx \subset \Real^{\nx}$, and $\f:\Dx \rightarrow \Real^{\nx}$. Suppose functional form of $\f$ is unknown.
However, we have been given ordered pairs $\{\y_i, \ydot_i\}$ for $i = 1,2\cdots \nd$, where $D:=\{\y_i\}_{i=1}^{\nd} \subset \Dx$.
Given data points in $D$ can come from a particular trajectory or random sparse observations. Our objective is to model unknown dynamics based upon available data.
In real world systems, often $\ydot_i$ is not available, and must be approximated from available $\y_i$. Approaches discussed in \cite{sindy} to address this issue may be adopted here as well.
However, for the purpose of this preliminary study, we assume that $\ydot_i$ are known.

Dynamics equation in its component form can be written as follows
\begin{align*}
  \xdot = \begin{pmatrix} \dx_1 \\ \dx_2 \\ \vdots \\ \dx_{\nx} \end{pmatrix}
        = \begin{pmatrix} f_1(\x) \\ f_2(\x) \\ \vdots \\ f_{\nx}(\x) \end{pmatrix}
        = \f(\x),
\end{align*}
where $f_i:\Dx \rightarrow \Real$ is a scalar function. Then we can use the formulation developed in Section \ref{fcn_approx} to construct a \maxent based approximate function $\fhat_i(\cdot)$ corresponding to unknown function $f_i(\cdot)$. Therefore, approximated dynamics model is given by
\begin{align*}
\xdot = \f(\x) \approx \hat{\f}(\x) := \begin{pmatrix} \fhat_1(\x) \\ \fhat_2(\x) \\ \vdots \\ \fhat_{\nx}(\x) \end{pmatrix}.
\end{align*}
We briefly summarize the framework discussed in Section \ref{maxent_fcn} and \ref{fcn_approx} in the context of dynamics modeling as the algorithm below.

\begin{algorithm}
\SetAlgoLined
 \tcp{ \small  Construct approximant from given data}
 Given: $\ydot_i$, for $\y_i \in D \subset \Dx$, $i = 1,2 \cdots \nd$. \\
 Select $\{\x_i\}_{i=1}^{\nb}=:B \subset \Dx$ such that $\Dx = \cvxh{B}$. \\
 \For{$i=1,2 \cdots \nd$}{
  Evaluate and store $\mb{\Psi}(\y_i)$ w.r.t. nodes in $B$ by solving \eqn{PhiME_loc}.\\
  }

  \For{$i=1,2 \cdots \nx$}{
   Determine and store $\q_i^*$ corresponding to scalar function $f_i(\cdot)$ by solving \eqn{opt_a}.\\
   }
 \tcp{\small Evaluate approximant for a given input}
    For any $\x \in \Dx$, evaluate $\mb{\Psi}(\x)$ w.r.t. nodes in $B$ by solving \eqn{PhiME_loc}.\\
    $\hat{\f}(\x) = \begin{pmatrix} \fhat_1(\x) \\ \fhat_2(\x) \\ \vdots \\ \fhat_{\nx}(\x) \end{pmatrix}
    = \begin{pmatrix} \q_1^{*T}\mb{\Psi}(\x) \\ \q_2^{*T}\mb{\Psi}(\x) \\ \vdots \\ \q_{\nx}^{*T}\mb{\Psi}(\x) \end{pmatrix}$
 \caption{Data driven modeling of dynamics using \maxent basis functions}
 \label{algo}
\end{algorithm}

\subsection*{Examples}
\subsubsection{Lorenz system} This is a system of three nonlinear ordinary differential equations given by \cite{lorenz}
\begin{align*}
\dx_1 &= \sigma (x_2-x_1),\\
\dx_2 &= x_1(\rho-x_3) - x_2,\\
\dx_3 &= x_1 x_2 - \gamma x_3,
\end{align*}
where, $\sigma = 10, \rho =28, \gamma = 8/3$. We assume that we have been given a trajectory with 500 observations, i.e. $\nd=500$. We select a uniform grid of $5^3$ points for basis nodes, such that all given observations lie within the convex hull of basis nodes.
Generally a given trajectory does not span the entire possible domain of the states. Therefore, for better accuracy of the approximated model for given data set, we select a subset of given data points as basis nodes. In this particular case we use sparsely distributed $100$ data points as basis nodes, therefore, $\nb = 5^3+100$.
Then Algorithm \ref{algo} is employed to construct an approximate model for the dynamics, which is integrated numerically to generate a trajectory. Its comparison with true trajectory is shown in \fig{lorenz}.
RMS error for the approximated trajectory calculated w.r.t. true trajectory is $\mathcal{O}(10^{-4})$. SINDy exactly recovers the dynamics as equations are $2^{\text{nd}}$ order polynomials.

\begin{figure}[thpb]
   \centering
   \includegraphics[width=0.5\textwidth]{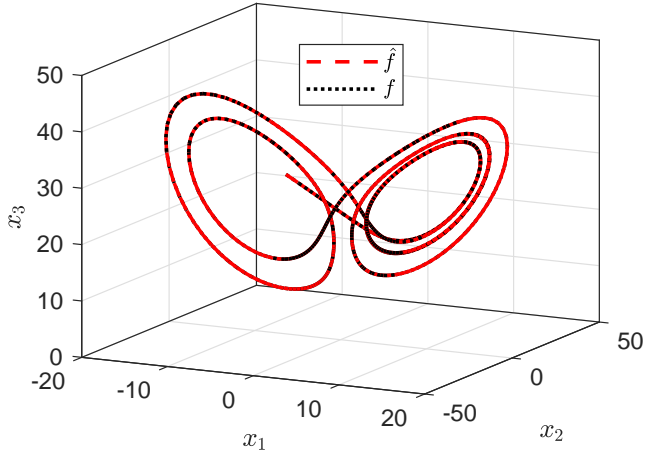}
   \caption{Comparison of true and approximate trajectories for Lorenz system.}
   \figlabel{lorenz}
\end{figure}

\subsubsection{Planar satellite orbit} We consider the dynamics of a satellite in a planar elliptic orbit around the Earth.
The state vector is defined by the radial distance of the satellite from Earth ($r$), rate of change of radial distance ($\dot{r}$), true anomaly ($\theta$), and rate of change of true anomaly ($\dot{\theta}$). Therefore, $\nx=4$.
Governing equations of motion are discussed in \cite{curtis}. Second order differential equations are given by
\begin{align*}
\ddot{r} &= \frac{\mu \epsilon}{h_0} \dot{\theta} \cos(\theta) , \\
\ddot{\theta} &= -\frac{2 \dot{\theta} \dot{r}}{r},
\end{align*}
where $\mu$ is standard gravitational parameter of the Earth, $\epsilon$ is eccentricity of the orbit, and $h_0$ is magnitude of angular momentum of the satellite which is a constant of motion.

We simulate this system for two cases, namely, $\nd = 500$ and $\nd = 20$. We follow similar approach as in the previous example. For $\nd = 500$, we select $\nb = 5^4 + 100$ basis nodes. Results obtained for this case using Algorithm \ref{algo} and SINDy are compared in \fig{orbit}.
Top plot in \fig{orbit} shows the normalized trajectory of satellite for two orbital periods ($T$) in cartesian coordinates, where $R_E$ is radius of the Earth. Clearly, \maxent based dynamics model gives quite accurate trajectory. On the other hand, we notice significant deviation from the true trajectory for SINDy based modeling.

Bottom plot in \fig{orbit} shows evolution of the error in angular momentum ($h$) with time in normalized quantities. The \maxent based model demonstrates better accuracy than SINDy, especially in the second orbital period.

For the second case of $\nd = 20$, when we have sparse observations, which is usually the case for space objects, we use $\nb = 5^4 + 20$. As shown in \fig{orbit_sparse}, SINDy based model deviates from the true trajectory within $10 \%$ of the first orbital period. On the other hand, we observe that \maxent based model predicts the trajectory quite accurately.
Although not shown in the figure, normalized error in angular momentum for \maxent based model is less than $10^{-2}$ for this case. This example demonstrates that \maxent based approach can model dynamics quite accurately even if the the available data is sparse.

\begin{figure}[h!]
   \centering
   \includegraphics[width=0.5\textwidth]{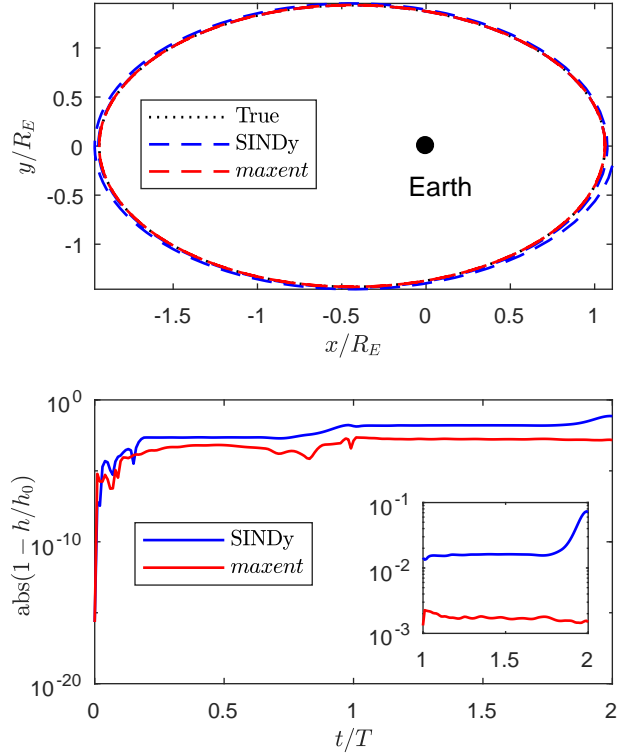}
   \caption{Comparison of true and approximate satellite trajectories, $\nd = 500$ and $\nb = 5^4+100$.}
   \figlabel{orbit}
\end{figure}

\begin{figure}[h!]
   \centering
   \includegraphics[width=0.5\textwidth]{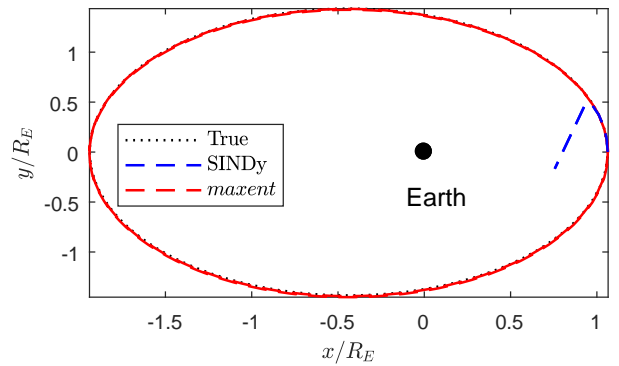}
   \caption{Comparison of true and approximate satellite trajectories using sparse observations, $\nd = 20$ and $\nb = 5^4+20$.}
   \figlabel{orbit_sparse}
\end{figure}

\section{CONCLUSIONS} \label{conclude}
In this paper we presented a framework to construct approximant from scattered data. We used maximum entropy functions which are derived from the data as bases for this purpose.

We demonstrated the application of our framework for few function approximation examples and dynamic systems. Numerical results show that \maxent based framework can be used for constructing approximants with reasonable accuracy, even if the available data set is sparse and no information is available about the functional form of underlying governing equations.

Our function approximation formulation has three control parameters, namely, $\nb$ which is the number of basis nodes and hence basis functions, $\beta$ which controls the degree of locality of basis functions, and $\alpha$ which regularizes the optimal basis coefficients $\q^*$.
As demonstrated in one of the examples, $\alpha$ can be used to make the basis set sparse, consequently reducing the computational complexity of the approximant.
However, how to determine optimal parameters that will result in the most accurate approximant is a topic of ongoing investigation.
Perhaps, established techniques for tuning parameters from machine learning literature can be adopted for this purpose, and will be addressed in our future works.

\addtolength{\textheight}{-12cm}   % This command serves to balance the column lengths

%%%%%%%%%%%%%%%%%%%%%%%%%%%%%%%%%%%%%%%%%%%%%%%%%%%%%%%%%%%%%%%%%%%%%%%%%%%%%%%%

%%%%%%%%%%%%%%%%%%%%%%%%%%%%%%%%%%%%%%%%%%%%%%%%%%%%%%%%%%%%%%%%%%%%%%%%%%%%%%%%


% \section*{ACKNOWLEDGMENT}

% Put sponsor acknowledgments in the unnumbered footnote on the first page.

%%%%%%%%%%%%%%%%%%%%%%%%%%%%%%%%%%%%%%%%%%%%%%%%%%%%%%%%%%%%%%%%%%%%%%%%%%%%%%%%

\bibliographystyle{./IEEEtran}
\bibliography{root}

\end{document}